\documentclass[a4paper,11pt,twoside]{article}
\usepackage{makeidx}
\usepackage{amsfonts}
\usepackage{fullpage}
\usepackage{amsmath}
\usepackage{amssymb}
\usepackage{amsthm}
\usepackage{appendix}
\usepackage{manfnt}
\usepackage{marginnote}
\usepackage{xfrac}
\newtheorem{thm}{Theorem}[section]
\newtheorem{lem}[thm]{Lemma}

\newtheorem{cor}[thm]{Corollary}
\theoremstyle{definition}
\newtheorem{ex}[thm]{Example}
\newtheorem{defn}[thm]{Definition}

\theoremstyle{remark}
\newtheorem{rmk}[thm]{Remark}
\usepackage{theoremref}
\usepackage[all]{xy}
\usepackage{mathrsfs}
\newcommand{\op}{\mathrm{op}}
\newcommand{\Pes}{\mathrm{Pes}}

\newcommand{\fp}{\mathrm{fp}}

\newcommand{\Coker}{\mathrm{Coker}}
\newcommand{\Ker}{\mathrm{Ker}}
\newcommand{\Img}{\mathrm{Im}}

\newcommand{\Ab}{\mathrm{Ab}}

\newcommand{\Ext}{\mathrm{Ext}}
\newcommand{\Tor}{\mathrm{Tor}}
\newcommand{\rmods}{\mathrm{mod}\text{-}}
\newcommand{\lmods}{\text{-}\mathrm{mod}}

\newcommand{\Lmods}{\text{-}\mathrm{Mod}}

\begin{document}
\title{The Auslander-Gruson-Jensen Recollement}
\author{Samuel Dean and Jeremy Russell}
\date{}
\maketitle
\abstract{For any ring $R$, the Auslander-Gruson-Jensen functor 
$$D_A:\fp(R\Lmods,\Ab)\to (\rmods R,\Ab)^\op$$
is the exact functor which sends a representable functor $(X,-)$ to the tensor functor $-\otimes X$. We show that this functor admits a fully faithful right adjoint $D_R$ and a fully faithful left adjoint $D_L$. That is, we show that $D_A$ is part of a recollement of abelian categories. In particular, this shows that $D_A$ is a localisation and a colocalisation which gives an equivalence of categories
\begin{displaymath}
\frac{\fp(R\Lmods,\Ab)}{\{F:D_AF=0\}}\simeq (\rmods R,\Ab)^\op.
\end{displaymath}
We show that $\{F:D_AF=0\}$ is the Serre subcategory of $\fp(R\Lmods,\Ab)$ consisting of finitely presented functors which arise from a pure-exact sequence. As an application of our main result, we show that the 0-th right pure-derived functor of a finitely presented functor $R\Lmods\to\Ab$ is also finitely presented.}
\section{Introduction}
All categories and functors mentioned in this paper are additive. Throughout the paper, let $R$ denote an arbitrary ring.

The category $R\Lmods$ is the category of left $R$-modules, and $\rmods R$ is the category of finitely presented right $R$-modules. The category $\fp(R\Lmods,\Ab)$ of finitely presented functors consists of all functors $F:R\Lmods\to\Ab$ such that there are left $R$-modules $X$ and $Y$ and there exists an exact sequence of functors
\begin{displaymath}
\xymatrix{(Y,-)\ar[r]& (X,-)\ar[r]&F\ar[r]& 0.}
\end{displaymath}
The category $(\rmods R,\Ab)$ is the category of all functors $\rmods R\to\Ab$.

Both $\fp(R\Lmods,\Ab)$ and $(\rmods R,\Ab)$ are abelian, with kernels and cokernels computed object-wise (see \cite[Section 2]{auslander1965} for this).

A \textbf{recollement of abelian categories} is a situation consisting of additive functors 
\begin{displaymath}
\xymatrix{\mathcal{A}'\ar[rr]|{i_*}&&\mathcal{A}\ar@<2.5ex>[ll]^{i^!}\ar@<-2.5ex>[ll]_{i^*}\ar[rr]|{j^*}&&\mathcal{A}''\ar@<2.5ex>[ll]^{j_*}\ar@<-2.5ex>[ll]_{j_!}}
\end{displaymath}
between abelian categories $\mathcal{A}'$, $\mathcal{A}$ and $\mathcal{A}''$ such that the following hold:
\begin{list}{$\bullet$}{•}
\item $\Img(i_*)=\Ker (j^*)$.
\item $i_*$ is fully faithful and $i^*\dashv i_*\dashv i^!$.
\item $j_!\dashv j^*\dashv j_*$ and $j_!$ and $j_*$ are fully faithful.
\end{list}
The goal of this paper is to introduce the \textbf{Auslander-Gruson-Jensen recollement}
\begin{displaymath}
\xymatrix{\Ker(D_A)\ar[rr]|{\subseteq}&&\fp(R\Lmods,\Ab)\ar@<2.5ex>[ll]\ar@<-2.5ex>[ll]\ar[rr]|{D_A}&&(\rmods R,\Ab)^\op.\ar@<2.5ex>[ll]^{D_R}\ar@<-2.5ex>[ll]_{D_L}}
\end{displaymath}
The functor $D_A$ sends a finitely presented functor $F\in\fp(R\Lmods,\Ab)$ to the functor $D_AF\in(\rmods R,\Ab)$ which is given by $(D_AF)M=(F,M\otimes-)$ for any $M\in\rmods R$. We will show that $\Ker(D_A)$ consists of those finitely presented functors $F:R\Lmods\to\Ab$ which are part of an exact sequence
\begin{displaymath}
\xymatrix{0\ar[r]& (Z,-)\ar[r]&(Y,-)\ar[r]& (X,-)\ar[r]& F\ar[r]& 0}
\end{displaymath}
where 
\begin{displaymath}
\xymatrix{0\ar[r]& X\ar[r]& Y\ar[r]&Z\ar[r]& 0}
\end{displaymath}
is a pure-exact sequence.

The basic properties of recollements are well-known. We refer to \cite{comprecoll} for these, and to \cite[Section 3.4]{borceux1} for background on fully faithful adjoints. In particular, we use the fact that, to determine a recollement, in the above notation, one need only define the functor $j^*$, prove the existence of a left and right adjoint $j_!$ and $j_*$, and show that one of those adjoints is fully faithful. 

The Auslander-Gruson-Jensen functor is an extension of the well-known Auslander-Gruson-Jensen equivalence
\begin{displaymath}
d:\fp(R\lmods,\Ab)\to (\fp(\rmods R,\Ab))^\op
\end{displaymath}
which was first written down explicitly in the 1980s, by Auslander in \cite{auslander1986} and by Gruson and Jensen in \cite{gj}. It sends a finitely presented functor $F:R\lmods\to\Ab$ to the functor $dF:\rmods R\to\Ab$ defined by $(dF)M=(F,M\otimes-)$ for any finitely presented right $R$-module $M$. It is also determined by the fact that it is exact and exchanges representables and finitely presented tensor functors. Its inverse, 
\begin{displaymath}
d:(\fp(\rmods R,\Ab))^\op\to\fp(R\lmods,\Ab)
\end{displaymath}
has a similar definition.

The Auslander-Gruson-Jensen equivalence has a syntactic characterisation in the model theory of modules. The discovery of this connection was made in a series of contributions:
\begin{enumerate}
\item Prest defines the notion of a dual $D\phi$ of a pp formula $\phi$. See \cite{prest1988} for details.
\item Herzog defines a category of pp-pairs, and shows that Prest's duality of pp formulas extends to an equivalence of categories. See \cite{herzog1993} for details.
\item Herzog's categories of pp-pairs are shown by Burke to be equivalent to categories of finitely presented functors on finitely presented modules, and the Auslander-Gruson-Jensen equivalence is shown to be equivalent to elementary duality of pp-pairs. See \cite{burke} for details.
\end{enumerate}

We note that the recollement which we construct here is an analogue of the recollement
\begin{displaymath}
\xymatrix{\Ker(\textbf{R})\ar[rr]&&\fp((R\Lmods)^\op,\Ab)\ar@<2.5ex>[ll]\ar@<-2.5ex>[ll]\ar[rr]|{\textbf{R}}&&((R\lmods)^\op,\Ab).\ar@<2.5ex>[ll]\ar@<-2.5ex>[ll]}
\end{displaymath}
constructed by Krause in \cite{krause1998}. Here $\textbf{R}$ sends a functor $G\in\fp((R\Lmods)^\op,\Ab)$ to its restriction $G|_{(R\lmods)^\op}$. Its kernel, $\Ker(\textbf{R})$, consists of those finitely presented functors $G:(R\Lmods)^\op\to\Ab$  which are part of an exact sequence
\begin{displaymath}
\xymatrix{0\ar[r]& (-,X)\ar[r]&  (-,Y)\ar[r]&  (-,Z)\ar[r]& G \ar[r]&  0}
\end{displaymath}
where 
\begin{displaymath}
\xymatrix{0\ar[r]&  X\ar[r]&  Y\ar[r]&  Z\ar[r]&  0}
\end{displaymath}
is a pure-exact sequence.

Studying purity by using contravariant functors (i.e.~the category $((\rmods R)^\op,\Ab)$) can be generalised to the setting of locally finitely presented categories, as discussed by Crawley-Boevey in \cite[Section 3]{crawleyboevey1994}, where it is shown that, for a locally finitely presented category $\mathcal{C}$, there is a fully faithful functor $\mathcal{C}\to((\fp\mathcal{C})^\op,\Ab)$, and that, under this embedding, pure-exact sequences in $\mathcal{C}$ correspond to exact sequences between flat functors. In particular, this shows that the class of pure-exact sequences form an exact structure on $\mathcal{C}$. However, it is well-known that, in module categories, pure-exact sequences have an equivalent formulation in terms of tensor products. That is, a short exact sequence
\begin{displaymath}
\xymatrix{0\ar[r]&X\ar[r]&Y\ar[r]&Z\ar[r]&0}
\end{displaymath}
is a pure-exact sequence in $R\Lmods$ if and only if 
\begin{displaymath}
\xymatrix{0\ar[r]&-\otimes X\ar[r]&-\otimes Y\ar[r]&-\otimes Z\ar[r]&0}
\end{displaymath}
is exact in $(\rmods R,\Ab)$. For this, see for example Crawley-Boevey's more general statement \cite[Section 3.3, Lemma 4]{crawleyboevey1994}. This begs the question of whether there is a recollement which corresponds to this alternative characterisation of purity for modules. This is what we provide.

Various interesting functors, although not necessarily finitely presented when restricted to $R\lmods$, are finitely presented as functors $R\Lmods\to\Ab$. Obvious examples include functors of the form $\Ext^n(A,-):R\Lmods\to\Ab$ where $A$ is an arbitrary left $R$-module, for any $n\geqslant 0$. Also, important functors can be lost by restricting to finitely presented modules. An illustrative example is the functor $(\mathbb{Q},-):\Ab\to\Ab$, which is in $\fp(\Ab,\Ab)$, but $(\mathbb{Q},-)|_{\rmods\mathbb{Z}}=0$. This is motivation for moving to $\fp(R\Lmods,\Ab)$.
\subsubsection*{Acknowledgements}Both authors would like to sincerely thank Rosie Laking for asking certain question which ultimately led to the consideration of $D_A$ and the discovery that it admitted both adjoints. They would also like to thank Henning Krause for pointing out that we had shown that $D_A$ is a Serre quotient functor.

The first author would like to thank his supervisor, Mike Prest, for his on-going support. He would also like to thank Lidia Angeleri H\"ugel and David Pauksztello for useful suggestions.

The second author would like to thank Mike Prest in particular for supporting multiple visits to the University of Manchester which was fundamental in discovering this problem. He would also like to thank Ivo Herzog for hosting him while visiting Ohio State University at Lima and making him aware of the connection between model theory and representation theory.

The authors would like to sincerely thank the referee for all of the suggestions made during the peer review process.  As a result of these suggestions, the authors have improved drastically the exposition of the main result and significantly reduced the length of the paper.  The original paper could not have been easy to read and we greatly appreciate the fact that the referee took the time to seriously consider how to better present the results.
\section{The Auslander-Gruson-Jensen functor}
In this section, we define the Auslander-Gruson-Jensen functor and prove that it is exact.
\begin{defn}The \textbf{Auslander-Gruson-Jensen functor} 
\begin{displaymath}
D_A:\fp(R\Lmods,\Ab)\to(\rmods R,\Ab)^\op
\end{displaymath}
is given by 
\begin{displaymath}
(D_AF)M=(F,M\otimes-)
\end{displaymath}
for any $F\in\fp(R\Lmods,\Ab)$ and any $M\in\rmods R$.
\end{defn}
It follows directly from the Yoneda lemma that, for any $X\in R\Lmods$, there is an isomorphism $D_A(X,-)\cong -\otimes X$ which is natural in $X$.
\begin{ex}\thlabel{example}For any left $R$-module $X$ and any positive integer $n$, there is an isomorphism 
$$D_A\Ext^n(X,-)\cong \Tor_n(-,X)$$
which is natural in $X$. To see this, we use Yoneda's result that, for any right exact functor $G$, there is an isomorphism
\begin{displaymath}
(\Ext^n(X,-),G)\cong (L_nG)X
\end{displaymath}
which is natural in $X$ and $G$ (see, \cite{yoneda}, \cite[Theorem 1.2]{hiltonrees} or \cite[Corollary 5.4]{auslander1965} for this). In particular we have the well-known result that, for \emph{any} right $R$-module $M$ there is an isomorphism
\begin{displaymath}
(\Ext^n(X,-),M\otimes -)\cong L_n(M\otimes -)(X)=\Tor_n(M,X)
\end{displaymath}
which is natural in $X$ and $M$. Restricting to the special case in which $M$ is finitely presented gives the desired result.
\end{ex}
\begin{cor}$D_A$ is exact.
\end{cor}
\begin{proof}It is well-known and easily proved that a functor $G:R\Lmods\to\Ab$ is right exact if and only if $\Ext^1(F,G)=0$ for any $F\in\fp(R\Lmods,\Ab)$. In particular, $\Ext^1(F,M\otimes -)=0$ for any right $R$-module $M$.

Let 
\begin{displaymath}
\xymatrix{0\ar[r]&F'\ar[r]&F\ar[r]&F''\ar[r]&0}
\end{displaymath}
be an exact sequence in $\fp(R\Lmods,\Ab)$. Then, for any right $R$-module $M$, the sequence
\begin{displaymath}
\xymatrix{0\ar[r]&(F'',M\otimes-)\ar[r]&(F,M\otimes-)\ar[r]&(F',M\otimes-)\ar[r]&\Ext^1(F,M\otimes-)}
\end{displaymath}
is exact. Since $\Ext^1(F,M\otimes-)=0$, restricting to the special case where $M$ is finitely presented gives the desired result.
\end{proof}
\begin{lem}For any $M\in\rmods R$, there is an isomorphism $D_A(M\otimes-)\cong (M,-)$ which is natural in $M$.
\end{lem}
\begin{proof}
It is well-known that the tensor embedding $\rmods R\to \fp(R\Lmods,\Ab)$ is fully faithful. Therefore, for any $N\in \rmods R$, there is an isomorphism
\begin{displaymath}
D_A(M\otimes -)(N)=(M\otimes-, N\otimes -)\cong (M,N)
\end{displaymath}
which is natural in $M$ and $N$.
\end{proof}
\section{The right adjoint of $D_A$}
\begin{defn}For any $F\in (\rmods R,\Ab)$, we define the functor $D_RF:R\Lmods\to\Ab$ by $(D_RF)M=(F,-\otimes M)$ for any $M\in R\Lmods$.
\end{defn}
\begin{lem}\thlabel{fpfun}$D_RF\in\fp(R\Lmods,\Ab)$ for any $F\in(\rmods R,\Ab)$. 
\end{lem}
\begin{proof}
First, note that, for any indexed family $\{F_i\}_{i\in I}$ of functors in $\fp(R\Lmods,\Ab)$, the functor $\prod_{i\in I}F_i$ defined by $(\prod_{i\in I}F_i)M=\prod_{i\in I}(F_iM)$ for any $M\in R\Lmods$, is finitely presented. To see this, observe that for any indexed family of modules $\{X_i\}_{i\in I}$, there is an isomorphism $\prod_{i\in I}(X_i,-)\cong(\bigoplus_{i\in I}X_i,-)$ of functors $R\Lmods\to\Ab$. Since products are exact in $\Ab$, if $F_i$ is given by an exact sequence
\begin{displaymath}
\xymatrix{(Y_i,-)\ar[r]&(X_i,-)\ar[r]&F_i\ar[r]& 0}
\end{displaymath}
then there is an exact sequence
\begin{displaymath}
\xymatrix{\prod_{i\in I}(Y_i,-)\ar[r]&\prod_{i\in I}(X_i,-)\ar[r]&\prod_{i\in I}F_i\ar[r]& 0}
\end{displaymath}
and hence an exact sequence
\begin{displaymath}
\xymatrix{\left(\bigoplus_{i\in I}Y_i,-\right)\ar[r]& \left(\bigoplus_{i\in I}X_i,-\right)\ar[r]&\prod_{i\in I}F_i\ar[r]& 0,}
\end{displaymath}
so $\prod_{i\in I} F_i\in\fp(R\Lmods,\Ab)$.

For any finitely presented right $R$-module $X$, $X\otimes-\in\fp(R\Lmods,\Ab)$ (see \cite[Lemma 6.1]{auslander1965} for this), and therefore $D_R(X,-)\simeq X\otimes-\in\fp(R\Lmods,\Ab)$.

For any $F\in(\rmods R,\Ab)$ there is an exact sequence
\begin{displaymath}
\xymatrix{\bigoplus_{j\in J}(Y_j,-)\ar[r]&\bigoplus_{i\in I}(X_i,-)\ar[r]& F\ar[r]& 0.}
\end{displaymath}
Therefore, there is also an exact sequence
\begin{displaymath}
\xymatrix{0\to D_RF\ar[r]& \prod_{i\in I}D_R(X_i,-)\ar[r]&\prod_{j\in J}D_R(Y_j,-)}
\end{displaymath}
and since $\fp(R\Lmods,\Ab)$ is closed under kernels, our above discussion shows that $D_RF\in\fp(R\Lmods,\Ab)$.
\end{proof}
\thref{fpfun} shows that we have a functor 
\begin{displaymath}
D_R:(\rmods R,\Ab)^\op\to \fp(R\Lmods,\Ab).
\end{displaymath}

We will continuously use the following well-known fact, and its dual.
\begin{thm}\thlabel{extendinglemma}Let $\mathcal{A}$ be an abelian category with enough projectives, and let $\mathcal{P}$ be the full subcategory of $\mathcal{A}$ of projective objects in $\mathcal{A}$.
\begin{enumerate}
\item Let $\mathcal{B}$ be an additive category with cokernels and let $f:\mathcal{P}\to\mathcal{B}$ be an additive functor. There is, up to isomorphism, a unique right exact functor $F:\mathcal{A}\to\mathcal{B}$ such that $F|_\mathcal{P}\cong f$.
\item Let $\mathcal{B}$ be an additive category with cokernels and let $F,G:\mathcal{A}\to\mathcal{B}$ be right exact functors. If $F|_\mathcal{P}\cong G|_\mathcal{P}$ then $F\cong G$.
\end{enumerate}
\end{thm}
\begin{proof}
The first part is identical to \cite[Lemma 2.4]{krause1998} and the second part follows immediately from the first.
\end{proof}
For convenience, we shall state the dual of \thref{extendinglemma}.
\begin{thm}\thlabel{dextendinglemma}Let $\mathcal{A}$ be an abelian category with enough injectives, and let $\mathcal{I}$ be the full subcategory of $\mathcal{A}$ of injective objects in $\mathcal{A}$.
\begin{enumerate}
\item Let $\mathcal{B}$ be an additive category with kernels and let $f:\mathcal{I}\to\mathcal{B}$ be an additive functor. There is, up to isomorphism, a unique left exact functor $F:\mathcal{A}\to\mathcal{B}$ such that $F|_\mathcal{I}\cong f$.
\item Let $\mathcal{B}$ be an additive category with kernels and let $F,G:\mathcal{A}\to\mathcal{B}$ be left exact functors. If $F|_\mathcal{I}\cong G|_\mathcal{I}$ then $F\cong G$.
\end{enumerate}
\end{thm}
\begin{thm}\thlabel{radj}There is an adjunction $D_A\dashv D_R$. That is, for functors $F\in\fp(R\Lmods,\Ab)$ and $G\in(\rmods R,\Ab)$, there is an isomorphism 
\begin{displaymath}
(F,D_RG)\cong (G,D_AF)
\end{displaymath}
which is  natural in $G$ and $F$.
\end{thm}
\begin{proof}
Since the functors $(-,D_RG),(G,D_A-):(\fp(R\Lmods,\Ab))^\op\to \Ab$ are left exact, \thref{extendinglemma} tells us that we need only check that they agree on projectives in $\fp(R\Lmods,\Ab)$, which are precisely the representable functors. Indeed, for any $X\in R\Lmods$ and $G\in(\rmods R,\Ab)$, there are isomorphisms
\begin{displaymath}
((X,-),D_RG)\cong (D_RG)X\cong (G,-\otimes X)\cong (G,D_A(X,-))
\end{displaymath}
which are natural in $X$ and $G$, as required. This completes the proof.
\end{proof}
\section{The left adjoint of $D_A$}
It is well-known that any Grothendieck category has injective hulls (see, for example, \cite[Chapter III, Theorem 3.2]{mitchell}). Therefore, $(\rmods R,\Ab)$ has injective hulls. It was shown by Gruson and Jensen in \cite{gj} that a functor $F\in(\rmods R,\Ab)$ is injective if and only if $F\cong -\otimes U$ for a pure-injective left $R$-module $U$. 

We will construct a left adjoint $D_L:(\rmods R,\Ab)^\op\to\fp(R\Lmods,\Ab)$ of $D_A$. Since it would have an exact right adjoint, $D_L$ sends projective objects to projective objects. That is, the injective functors in $(\rmods R,\Ab)$ should be sent to representable functors in $\fp(R\Lmods,\Ab)$ by $D_L$. The most straightforward and obvious definition is to define $D_L$ to be the unique right exact functor which sends an injective functor $-\otimes N$ to $(N,-)$ for any pure-injective $N\in R\Lmods$. That is, on injective objects $D_L$ agrees with $Y\circ \text{ev}^\op_R$ where $\text{ev}_R:(\rmods R,\Ab)\to R\Lmods$ is the functor which sends a functor $F\in(\rmods R,\Ab)$ to $F(R_R)$ (which has the natural structure of a left $R$-module since it is a left $\text{End}(R_R)$-module and $\text{End}(R_R)\cong R$) and $Y:(R\Lmods)^\op\to\fp(R\Lmods,\Ab)$ is the Yoneda embedding.
\begin{defn}Define the functor $D_L:(\rmods R,\Ab)^\op\to\fp(R\Lmods,\Ab)$ by
$$D_L=L_0(Y\circ \text{ev}_R^\op)$$
where $\text{ev}_R:(\rmods R,\Ab)\to R\Lmods$ is evaluation at the ring and $Y:(R\Lmods)^\op\to\fp(R\Lmods,\Ab)$ is the Yoneda embedding.
\end{defn}
\begin{lem}The functor $D_L$ is fully faithful.
\end{lem}
\begin{proof}
Let $\alpha:D_LG\to D_LF$ be a natural transformation, for some $F,G\in(\rmods R,\Ab)$. 

There are pure-injectives $U,V,X,Y\in R\Lmods$ and exact sequences
\begin{displaymath}
\xymatrix{0\ar[r]&F\ar[r]&-\otimes U\ar[r]&-\otimes V}
\end{displaymath}
and 
\begin{displaymath}
\xymatrix{0\ar[r]&G\ar[r]&-\otimes X\ar[r]&-\otimes Y.}
\end{displaymath}
Since $D_L$ is right exact, there is an induced commutative diagram
\begin{displaymath}
\xymatrix{(Y,-)\ar[d]\ar[r]&(X,-)\ar[d]\ar[r]&D_LG\ar[r]\ar[d]^\alpha&0\\
(V,-)\ar[r]&(U,-)\ar[r]&D_LF\ar[r]&0}
\end{displaymath}
with exact rows, since $\alpha$ lifts to a morphism between projective presentations. By the Yoneda embedding, this morphism is given by a unique commutative square
\begin{displaymath}
\xymatrix{U\ar[d]\ar[r]&V\ar[d]\\X\ar[r]&Y.}
\end{displaymath}

This commutative square induces a commutative diagram 
\begin{displaymath}
\xymatrix{0\ar[r]&F\ar[d]_\beta\ar[r]&-\otimes U\ar[d]\ar[r]&-\otimes V\ar[d]\\
0\ar[r]&G\ar[r]&-\otimes X\ar[r]&-\otimes Y}
\end{displaymath}
for some morphism $\beta:F\to G$. Since $D_L$ is right exact, a simple diagram chase shows that $D_L\beta=\alpha$.

To show that $D_L$ is faithful, assume $\alpha=0$. We need to show that $\beta=0$. Since $\alpha=0$, the image of the morphism $(X,-)\to(U,-)$ is contained in the image of $(V,-)\to(U,-)$, and it follows by the fact that $(X,-)$ is projective that there is a morphism $(X,-)\to (V,-)$ such that the diagram
\begin{displaymath}
\xymatrix{&(X,-)\ar[dl]\ar[d]\\
(V,-)\ar[r]&(U,-)}
\end{displaymath}
commutes. This morphism is induced by a morphism $V\to X$ such that 
\begin{displaymath}
\xymatrix{U\ar[d]\ar[r]&V\ar[dl]\\X&}
\end{displaymath}
commutes. It follows that the composition of $\beta$ and the monomorphism $G\to-\otimes X$ is zero, and therefore $\beta=0$. This completes the proof.
\end{proof}
\begin{thm}\thlabel{ladj}$D_L\dashv D_A$. That is, for functors $F\in\fp(R\Lmods,\Ab)$ and $G\in (\rmods R,\Ab)$, there is an isomorphism $(D_LG,F)\cong(D_AF,G)$ which is natural in $F$ and $G$.
\end{thm}
\begin{proof}
Let 
\begin{displaymath}
\xymatrix{(Y,-)\ar[r]&(X,-)\ar[r]&F\ar[r]&0}
\end{displaymath}
be an exact sequence in $\fp(R\Lmods,\Ab)$ for some $X,Y\in R\Lmods$ and $F\in\fp(R\Lmods,\Ab)$. Since $D_A$ is exact, we obtain an exact sequence
\begin{displaymath}
\xymatrix{0\ar[r]&D_AF\ar[r]&-\otimes X\ar[r]&-\otimes Y.}
\end{displaymath}
Let $U\in R\Lmods$ be a pure-injective. The functor $-\otimes U$ is injective in $(\rmods R,\Ab)$, so by applying $(-,-\otimes U)$ and using the fact that the tensor embedding $t:R\Lmods\to(\rmods R,\Ab)$ is fully faithful, we obtain an exact sequence
\begin{displaymath}
\xymatrix{(Y,U)\ar[r]&(X,U)\ar[r]&(D_AF,-\otimes U)\ar[r]&0.}
\end{displaymath}
Therefore, by the Yoneda lemma and the definition of $D_L$, there are isomorphisms
$$(D_AF,-\otimes U)\cong FU\cong ((U,-),F)\cong (D_L(-\otimes U),F).$$
A simple diagram chase shows that these isomorphisms are natural in $F$ and $U$.

The functors $(D_AF,-),(D_L-,F):(\rmods R,\Ab)\to\Ab$ are left exact. Since these functors agree on injectives, it follows from \thref{dextendinglemma} that for any $G\in (\rmods R,\Ab)$ there is an isomorphism $(D_AF,G)\cong (D_LG,F)$ which is natural in $G$ and $F$.
\end{proof}
\section{The main result and some of its consequences}
\begin{thm}[Auslander-Gruson-Jensen recollement]\thlabel{agjrec}There is a recollement of abelian categories
\begin{displaymath}
\xymatrix{\Ker(D_A)\ar[rr]|{\subseteq}&&\fp(R\Lmods,\Ab)\ar@<2.5ex>[ll]\ar@<-2.5ex>[ll]\ar[rr]|{D_A}&&(\rmods R,\Ab)^\op.\ar@<2.5ex>[ll]^{D_R}\ar@<-2.5ex>[ll]_{D_L}}
\end{displaymath}
\end{thm}
\begin{proof}
We have already constructed the right hand side of the diagram. The left hand side is uniquely determined by the right hand side (see \cite[Remark 2.8]{psavit} for details).
\end{proof}
\begin{cor}There is an equivalence of categories
\begin{displaymath}
\frac{\fp(R\Lmods,\Ab)}{\{F:D_AF=0\}}\simeq(\rmods R,\Ab)^\op.
\end{displaymath}
\end{cor}
\begin{proof}
\thref{radj} and \thref{ladj} show that $D_A$ is simultaneously a localisation and a colocalisation. Therefore it is the Serre quotient 
\begin{displaymath}
\fp(R\Lmods,\Ab)\to \fp(R\Lmods,\Ab)/{\{F:D_AF=0\}}
\end{displaymath}
and so we have the desired equivalence.
\end{proof}
\begin{lem}\thlabel{kerda}For a functor $F\in\fp(R\Lmods,\Ab)$, the following are equivalent.
\begin{enumerate}
\item $D_AF=0$.
\item If 
\begin{displaymath}
\xymatrix{0\ar[r]&(Z,-)\ar[r]&(Y,-)\ar[r]&(X,-)\ar[r]&F\ar[r]&0}
\end{displaymath}
is an exact sequence then
\begin{displaymath}
\xymatrix{0\ar[r]&X\ar[r]&Y\ar[r]&Z\ar[r]&0}
\end{displaymath}
is a pure-exact sequence.
\item There is a pure-exact sequence
\begin{displaymath}
\xymatrix{0\ar[r]&X\ar[r]&Y\ar[r]&Z\ar[r]&0}
\end{displaymath}
and an exact sequence
\begin{displaymath}
\xymatrix{0\ar[r]&(Z,-)\ar[r]&(Y,-)\ar[r]&(X,-)\ar[r]&F\ar[r]&0}
\end{displaymath}.
\item For any pure-injective $N\in R\Lmods$, $FN=0$.
\end{enumerate}
\end{lem}
\begin{proof}
\begin{list}{*}{}
\item[1 implies 2:] Since $D_A$ is exact, there is an exact sequence
\begin{displaymath}
\xymatrix{0\ar[r]&D_AF\ar[r]&-\otimes X\ar[r]&-\otimes Y\ar[r]&-\otimes Z\ar[r]&0.}
\end{displaymath}
Since $D_AF=0$, this shows that $0\to X\to Y\to Z\to 0$ is pure-exact.

\item[2 implies 3:] Obvious since $\fp(R\Lmods,\Ab)$ has global dimension at most 2.

\item[3 implies 4:] Since $N$ is pure-injective, $(Y,N)\to (X,N)$ is an epimorphism, so
$$FN\cong \Coker((Y,N)\to (X,N))=0.$$

\item[4 implies 1:] For any pure-injective $N\in R\Lmods$, we have 
\begin{displaymath}
(D_AF,-\otimes N)\cong (D_L(-\otimes N),F)\cong ((N,-),F)\cong FN.
\end{displaymath}
Therefore $D_AF$ has no non-zero maps into an injective functor, so since $(\rmods R,\Ab)$ has enough injectives, $D_AF=0$.
\end{list}
\end{proof}
\begin{cor}There is an equivalence of categories
\begin{displaymath}
(\Ker(D_A))^\op\simeq \Pes(R\Lmods)
\end{displaymath}
where $\Pes(R\Lmods)$ is the category of pure-exact sequences and homotopy classes of chain maps.
\end{cor}
\begin{proof}
Let $\mathcal{P}$ be the category of pure-exact sequences and chain maps. There is a functor $\mathcal{P}\to(\Ker(D_A))^\op$ which sends a pure-exact sequence 
\begin{displaymath}
\xymatrix{0\ar[r]&X\ar[r]&Y\ar[r]&Z\ar[r]&0}
\end{displaymath}
to the functor $F$ which is defined by the exact sequence
\begin{displaymath}
\xymatrix{0\ar[r]&(Z,-)\ar[r]&(Y,-)\ar[r]&(X,-)\ar[r]&F\ar[r]&0.}
\end{displaymath}
Natural transformations between finitely presented functors in $\Ker(D_A)$ lift uniquely to homotopy classes of chain maps between projective resolutions, and hence to homotopy classes of chain maps between pure-exact sequences. Therefore, the chain maps which are sent to zero by the functor $\mathcal{P}\to(\Ker(D_A))^\op$ are those in the ideal $\mathcal{I}$ of null homotopic chain maps. By \thref{kerda}, $\mathcal{P}\to(\Ker(D_A))^\op$ is essentially surjective so
$$(\Ker(D_A))^\op\simeq \mathcal{P}/\mathcal{I}=\Pes(R\Lmods).$$\end{proof}
\begin{ex}As an example of how our results can be used, we shall discuss 0-th derived functors with respect to the pure-exact structure on $R\Lmods$. See \cite{buhler} for the theory of derived functors with respect to an exact structure. 

Consider the functor
$$\widetilde{L}_0Y:(R\Lmods)^\op\to\fp(R\Lmods,\Ab),$$
which is the 0-th left pure-derived functor of the Yoneda embedding 
$$Y:(R\Lmods)^\op\to\fp(R\Lmods,\Ab).$$ That is, $\widetilde{L}_0Y$ is the 0-th left derived functor of $Y$ taken with respect to the pure-exact structure on $R\Lmods$. Let $M$ be any left $R$-module, and find a pure-injective copresentation
\begin{displaymath}
\xymatrix{0\ar[r]&M\ar[r]&U\ar[r]&V}
\end{displaymath}
of $M$ (that is, an injective copresentation of $M$ with respect to the pure-exact structure on $R\Lmods$). This induces an injective copresentation
\begin{displaymath}
\xymatrix{0\ar[r]&-\otimes M\ar[r]&-\otimes U\ar[r]&-\otimes V}
\end{displaymath}
of $-\otimes M$ in $(\rmods R,\Ab)$. Since $D_L$ is right exact, there is an exact sequence
\begin{displaymath}
\xymatrix{(V,-)\ar[r]& (U,-)\ar[r]& D_L(-\otimes M)\ar[r]&0.}
\end{displaymath}
Therefore, $(\widetilde{L}_0Y)M=D_L(-\otimes M)$. It follows that the functor $\widetilde{L}_0Y$ is fully faithful, since it is the composition of the fully faithful embedding 
$$R\Lmods\to(\rmods R,\Ab):M\mapsto -\otimes M$$
and $D_L$. 

For any functor $F:R\Lmods\to \Ab$, the sequence
\begin{displaymath}
\xymatrix{0\ar[r]&((\widetilde{L}_0Y)M,F)\ar[r]&((U,-),F)\ar[r]&((V,-),F)}
\end{displaymath}
is exact. However, the sequence
\begin{displaymath}
\xymatrix{0\ar[r]&(\widetilde{R}^0F)M\ar[r]&FU\ar[r]&FV}
\end{displaymath}
is also exact by definition, where $\widetilde{R}^0F$ is the 0-th right pure-derived functor of $F$. By the Yoneda lemma and the exactness of these sequences, we have $(\widetilde{R}^0F)M=((\widetilde{L}_0Y)M,F)$.

For a functor $F\in\fp(R\Lmods,\Ab)$ it follows from the above discussion that, for any left $R$-module $M$, we have an isomorphism
\begin{align*}
(\widetilde{R}^0F)M&=((\widetilde{L}_0Y)M,F)
\\&=(D_L(-\otimes M),F)
\\&\cong (D_AF,-\otimes M)
\\&=(D_RD_AF)M
\end{align*}
which is natural in $M$. Therefore, $\widetilde{R}^0F\cong D_RD_AF$. There are two immediate consequences of this result. The first is that, for any $F\in\fp(R\Lmods,\Ab)$, $\widetilde{R}^0F\in\fp(R\Lmods,\Ab)$. The second is that, for any $F\in\fp(R\Lmods,\Ab)$, $D_RD_AF\cong F$ if and only if $F$ is left pure-exact, meaning that it send pure-exact sequences to left exact sequences. This is because $\widetilde{R}^0F$ is the unique left pure-exact functor which agrees with $F$ on pure-injectives, so $\widetilde{R}^0F\cong F$ if and only if $F$ is left pure-exact.

Let $\mathcal{L}$ be the full subcategory of $\fp(R\Lmods,\Ab)$ consisting of left pure-exact finitely presented functors. For any $G\in(\rmods R,\Ab)$, $G\cong D_AD_RG$ (since $D_R$ is fully faithful, the counit of adjunction $D_AD_R\to 1$ is an isomorphism by \cite[IV.3, Theorem 1]{maclane1998}) and therefore
\begin{displaymath}
D_RG\cong D_RD_AD_RG\cong \widetilde{R}^0(D_RG),
\end{displaymath}
so $D_RG\in\mathcal{L}$. Therefore, since $D_R$ is fully faithful and maps every functor in $(\rmods R,\Ab)$ into $\mathcal{L}$, and since, by the above analysis, any left pure-exact finitely presented functor is isomorphic to a functor in the image of $D_R$, the functor $D_R$ restricts to an equivalence $(\rmods R,\Ab)^\op\simeq \mathcal{L}$.
\end{ex}
\begin{rmk}In the introduction, we mentioned Krause's recollement, which is analogous to the Auslander-Gruson-Jensen recollement. This raises an interesting question. What is the relationship between these recollements? We will not answer this question here, but it will be addressed in \cite{dean3} where we show that these two recollements are instances of a general construction which uses 0-th derived functors. Using results from \cite{stabmodth}, we will actually show that these recollements can be restricted (in a non-trivial way) until the functors $\textbf{R}$ and $D_A$ become dual to each other, where $\textbf{R}:\fp((R\Lmods)^\op,\Ab)\to((R\lmods)^\op,\Ab)$ is the restriction functor.

In \cite{dean3}, we will also show that other recollements arise from the same methods, including recollements which correspond to the functors
$$w:\fp(R\Lmods,\Ab)\to (R\Lmods)^\op$$
and 
$$\text{ev}_R:(\rmods R,\Ab)\to R\Lmods.$$
\end{rmk}
\begin{rmk}Because $D_R$ and $D_L$ agree on injectives and $D_L$ is right exact, we have $D_L\cong L_0D_R$. Therefore, there is a canonical morphism $D_L\to D_R$. We can form an exact sequence
\begin{displaymath}
\xymatrix{0\ar[r]& K\ar[r]& D_L\ar[r]& D_R\ar[r]& \underline{D_R}\ar[r]& 0.}
\end{displaymath}
This is an example of the projective-stabilisation sequence introduced by Auslander and Bridger in \cite{stabmodth}. Note that we could have constructed this sequence from general results about recollements (see \cite[Proposition 4.4]{comprecoll}).
\end{rmk}

\end{document}